\documentclass{article}
\usepackage{amssymb,latexsym}
\usepackage{amsmath}
\usepackage{graphics}
\usepackage{graphicx}
\newtheorem{theorem}{Theorem}

\newtheorem{definition}[theorem]{Definition}

\numberwithin{equation}{section} \numberwithin{theorem}{section}
\numberwithin{corollary}{section} \numberwithin{lemma}{section}
\numberwithin{remark}{section} \numberwithin{notation}{section}

\begin{document}
\title{Smarandache Curves According to Sabban Frame on ${S^2}$}
\author{Kemal Ta\c{s}k\"{o}pr\"{u}, Murat Tosun \\
\small Faculty of Arts and Sciences, Department of Mathematics,\\
\small Sakarya University, Sakarya, 54187, TURKEY}

\date{}

\maketitle

\noindent \textbf{Abstract:} In this paper, we introduce special Smarandache curves according
to Sabban frame on ${S^2}$ and we give some characterization of Smarandache
curves. Besides, we illustrate examples of our results. \\


\section{Introduction}\label{S:intro}
The differential geometry of curves is usual starting point of students in field of differential geometry which is the field concerned with studying curves, surfaces, etc. with the use the concept of derivatives in calculus. Thus, implicit in the discussion, we assume that the defining functions are sufficiently differentiable, i.e., they have no concerns of cusps, etc. Curves are usually studied as subsets of an ambient space with a notion of equivalence. For example, one may study curves in the plane, the usual three dimensional space, curves on a sphere, etc. There are many important consequences and properties of curves in differential geometry. In the light of the existing studies, authors introduced new curves. Special Smarandache curves are one of them. A regular curve in Minkowski space-time, whose position vector is composed by Frenet frame vectors on another regular curve, is called a Smarandache curve \cite{Tmy5}. Special Smarandache curves have been studied by some authors \cite{Aat2,Çtk3,Byc4,Tmy5,Nur6}.\\\\
In this paper, we study special Smarandache curves such as  $\gamma t$, $td$, $\gamma td$ - Smarandache curves according to Sabban frame in Euclidean unit sphere ${S^2}$. We hope these results will be helpful to mathematicians who are specialized on mathematical modeling.

\section{Preliminaries}\label{S:intro}
The Euclidean 3-space $E^3$  provided with the standard flat metric given by
\[\left\langle , \right\rangle  = dx_1^2 + dx_2^2 + dx_3^2\]
where $(x_{1},x_{2},x_{3})$ is a rectangular coordinate system of $E^3$. Recall that, the norm of an arbitrary vector $X\in E^3$ is given by $\|X\|=\sqrt{\langle X,X \rangle}$. The curve $\alpha$ is called a unit speed curve if velocity vector  $\alpha'$ of $\alpha$  satisfies $\|\alpha'\|=1$. For vectors $v,w\in E^3$, it is said to be orthogonal if and only if $\left\langle v,w \right\rangle=0$. The sphere of radius $r=1$ and with center in the origin in the space $E^3$ is defined by
\[\begin{array}{l}S^2=\{P=(P_{1},P_{2},P_{3})|\langle P,P \rangle=1\}.\end{array}\]

\noindent Denote by $\left\{{T,N,B}\right\}$  the moving Frenet frame along the curve $\alpha$ in $E^3$. For an arbitrary curve $\alpha\in E^3$ , with first and second curvature, $\kappa$  and  $\tau$ respectively, the Frenet formulae is given by  \cite{Dcm1}
\[\begin{array}{l}
T' = \kappa N\\
N' =  - \kappa T + \tau B\\
B' =  - \tau N.
\end{array}\]

\noindent Now, we give a new frame different from Frenet frame. Let $\gamma$ be a unit speed spherical curve. We denote $s$ as the arc-length parameter of $\gamma$. Let us denote $t\left( s \right) = \gamma'\left( s \right)$, and we call $t\left( s \right)$ a unit tangent vector of $\gamma$. We now set a vector $d\left( s \right) = \gamma \left( s \right) \wedge t\left( s \right)$ along $\gamma$. This frame is called the Sabban frame of $\gamma$ on ${S^2}$ (Sphere of unit radius). Then we have the following spherical Frenet formulae of $\gamma$ :
\[\begin{array}{l}
\gamma ' = t\\
t' =  - \gamma  + {\kappa _g}d\\
d' =  - {\kappa _g}t
\end{array}\]
where   is called the geodesic curvature of $\kappa_{g}$  on $S^2$ and $\kappa_{g}=\langle t',d \rangle$ \cite{Koe7}.

\section{Smarandache Curves According to Sabban Frame on $S^2$}\label{S:intro}
In this section, we investigate Smarandache curves according to the Sabban frame on $S^2$ .Let $\gamma=\gamma(s)$ and $\beta=\beta(s^{*})$  be a unit speed regular spherical curves on $S^2$, and  $\{\gamma,t,d\}$ and $\{\gamma_{\beta},t_{\beta},d_{\beta}\}$ be the Sabban frame of these curves, respectively.

\subsection{$\gamma t$-Smarandache Curves}

\begin{definition}
Let $S^2$ be a unit sphere in $E^3$ and suppose that the unit speed regular curve $\gamma=\gamma(s)$ lying fully on $S^2$. In this case, $\gamma t$ - Smarandache curve can be defined by

\begin{equation}\label{2.1}
\beta(s^{*})=\frac{1}{\sqrt{2}}(\gamma+t).
\end{equation}
\end{definition}
Now we can compute Sabban invariants of $\gamma t$ - Smarandache curves. Differentiating the equation (3.1) with respect to $s$, we have
\[\beta '\left( {{s^ * }} \right) = \frac{{d\beta }}{{d{s^ * }}}\frac{{d{s^ * }}}{{ds}} = \frac{1}{{\sqrt 2 }}\left( {\gamma ' + t'} \right)\]
and

\[{t_\beta }\frac{{d{s^ * }}}{{ds}} = \frac{1}{{\sqrt 2 }}\left( {t + {\kappa _g}d - \gamma } \right)\]
where

\begin{equation}\label{2.1}
\frac{{d{s^ * }}}{{ds}} = \sqrt {\frac{{2 + {\kappa _g}^2}}{2}}.
\end{equation}
Thus, the tangent vector of curve $\beta$ is to be

\begin{equation}\label{2.1}
{t_\beta } = \frac{1}{{\sqrt {2 + {\kappa _g}^2} }}\left( { - \gamma  + t + {\kappa _g}d} \right).
\end{equation}
Differentiating the equation (3.3) with respect to $s$, we get

\begin{equation}\label{2.1}
{t_\beta }^\prime \frac{{d{s^ * }}}{{ds}} = \frac{1}{{{{\left( {2 + {\kappa _g}^2} \right)}^{\frac{3}{2}}}}}\left( {{\lambda _1}\gamma  + {\lambda _2}t + {\lambda _3}d} \right)
\end{equation}
where

\[\begin{array}{l}
{\lambda _1} = {\kappa _g}{\kappa _g}^\prime  - {\kappa _g}^2 - 2\\
{\lambda _2} =  - {\kappa _g}{\kappa _g}^\prime  - 2 - 2{\kappa _g}^2 - {\kappa _g}^4\\
{\lambda _3} = 2{\kappa _g} + 2{\kappa _g}^\prime  + {\kappa _g}^3.
\end{array}\]
Substituting the equation (3.2) into equation (3.4), we reach

\begin{equation}\label{2.1}
{t_\beta }^\prime  = \frac{{\sqrt 2 }}{{{{\left( {2 + {\kappa _g}^2} \right)}^2}}}\left( {{\lambda _1}\gamma  + {\lambda _2}t + {\lambda _3}d} \right).
\end{equation}
Considering the equations (3.1) and (3.3), it easily seen that

\begin{equation}\label{2.1}
{d_\beta } = \beta  \wedge {t_\beta } = \frac{1}{{\sqrt {4 + 2{\kappa _g}^2} }}\left( {{\kappa _g}\gamma  + \left( { - 1 - {\kappa _g}} \right)t + 2d} \right).
\end{equation}
From the equation (3.5) and (3.6), the geodesic curvature of $\beta(s^{*})$ is

\[\begin{array}{l}
{\kappa _g}^\beta  = \left\langle {{t_\beta }^\prime ,{d_\beta }} \right\rangle \\\\
\,\,\,\,\,\,\,\,\,\,\,=\frac{1}{{{{\left( {2 + {\kappa _g}^2} \right)}^{\frac{3}{2}}}}}\left( {{\lambda _1}{\kappa _g} + {\lambda _2}\left( { - 1 - {\kappa _g}} \right) + 2{\lambda _3}} \right).
\end{array}\]

\subsection{$td$-Smarandache Curves}

\begin{definition}
Let $S^2$ be a unit sphere in $E^3$ and suppose that the unit speed regular curve $\gamma=\gamma(s)$ lying fully on $S^2$. In this case, $td$ - Smarandache curve can be defined by

\begin{equation}\label{2.1}
\beta(s^{*})=\frac{1}{\sqrt{2}}(t+d).
\end{equation}
\end{definition}
Now we can compute Sabban invariants of $td$ - Smarandache curves. Differentiating the equation (3.7) with respect to $s$, we have

\[\beta '\left( {{s^ * }} \right) = \frac{{d\beta }}{{d{s^ * }}}\frac{{d{s^ * }}}{{ds}} = \frac{1}{{\sqrt 2 }}\left( {t' + d'} \right)\]
and

\[{t_\beta }\frac{{d{s^ * }}}{{ds}} = \frac{1}{{\sqrt 2 }}\left( { - \gamma  + {\kappa _g}d - {\kappa _g}t} \right)\]
where

\begin{equation}\label{2.1}
\frac{{d{s^ * }}}{{ds}} = \sqrt {\frac{{1 + 2{\kappa _g}^2}}{2}}.
\end{equation}
In that case, the tangent vector of curve $\beta$ is as follows

\begin{equation}\label{2.1}
{t_\beta } = \frac{1}{{\sqrt {1 + 2{\kappa _g}^2} }}\left( { - \gamma  - {\kappa _g}t + {\kappa _g}d} \right).
\end{equation}
Differentiating the equation (3.9) with respect to $s$, it is obtained that

\begin{equation}\label{2.1}
{t_\beta }^\prime \frac{{d{s^ * }}}{{ds}} = \frac{1}{{{{\left( {1 + 2{\kappa _g}^2} \right)}^{\frac{3}{2}}}}}\left( {{\lambda _1}\gamma  + {\lambda _2}t + {\lambda _3}d} \right)
\end{equation}
where

\[\begin{array}{l}
{\lambda _1} = 2{\kappa _g}{\kappa _g}^\prime  + {\kappa _g} + 2{\kappa _g}^3\\
{\lambda _2} =  - 1 - {\kappa _g}^\prime  - 3{\kappa _g}^2 - 2{\kappa _g}^4\\
{\lambda _3} =  - {\kappa _g}^2 + {\kappa _g}^\prime  - 2{\kappa _g}^4.
\end{array}\]
Substituting the equation (3.8) into equation (3.10), we get

\begin{equation}\label{2.1}
{t_\beta }^\prime  = \frac{{\sqrt 2 }}{{{{\left( {1 + 2{\kappa _g}^2} \right)}^2}}}\left( {{\lambda _1}\gamma  + {\lambda _2}t + {\lambda _3}d} \right)
\end{equation}
Using the equations (3.7) and (3.9), we easily find

\begin{equation}\label{2.1}
{d_\beta } = \beta  \wedge {t_\beta } = \frac{1}{{\sqrt {2 + 4{\kappa _g}^2} }}\left( {{\kappa _g}\gamma  - t + \left( {1 + {\kappa _g}} \right)d} \right).
\end{equation}
So, the geodesic curvature of $\beta(s^{*})$  is as follows

\[\begin{array}{l}
{\kappa _g}^\beta  = \left\langle {{t_\beta }^\prime ,{d_\beta }} \right\rangle \\\\
\,\,\,\,\,\,\,\,\,\,\,=\frac{1}{{{{\left( {1 + 2{\kappa _g}^2} \right)}^{\frac{3}{2}}}}}\left( {{\lambda _1}{\kappa _g} - {\lambda _2} + {\lambda _3}\left( {1 + {\kappa _g}} \right)} \right).
\end{array}\]

\subsection{$\gamma td$-Smarandache Curves}

\begin{definition}
Let $S^2$ be a unit sphere in $E^3$ and suppose that the unit speed regular curve $\gamma=\gamma(s)$ lying fully on $S^2$. Denote the Sabban frame of $\gamma(s)$, $\left\{ {\gamma ,t,d} \right\}$. In this case,  $\gamma td$ - Smarandache curve can be defined by

\begin{equation}\label{2.1}
\beta(s^{*})=\frac{1}{\sqrt{2}}(\gamma+t+d).
\end{equation}
\end{definition}
Lastly, let us calculate Sabban invariants of  $\gamma td$ - Smarandache curves. Differentiating the equation (3.13) with respect to $s$, we have

\[\beta '\left( {{s^ * }} \right) = \frac{{d\beta }}{{d{s^ * }}}\frac{{d{s^ * }}}{{ds}} = \frac{1}{{\sqrt 3 }}\left( {\gamma ' + t' + d'} \right)\]
and

\[{t_\beta }\frac{{d{s^ * }}}{{ds}} = \frac{1}{{\sqrt 3 }}\left( {t - \gamma  + {\kappa _g}d - {\kappa _g}t} \right)\]
where

\begin{equation}\label{2.1}
\frac{{d{s^ * }}}{{ds}} = \sqrt {\frac{{2\left( {1 - {\kappa _g} + {\kappa _g}^2} \right)}}{3}}.
\end{equation}
Thus, the tangent vector of curve $\beta$ is

\begin{equation}\label{2.1}
{t_\beta } = \frac{1}{{\sqrt {2\left( {1 - {\kappa _g} + {\kappa _g}^2} \right)} }}\left( { - \gamma  + \left( {1 - {\kappa _g}} \right)t + {\kappa _g}d} \right).
\end{equation}
Differentiating the equation (3.15) with respect to $s$, it is obtained that

\begin{equation}\label{2.1}
{t_\beta }^\prime \frac{{d{s^ * }}}{{ds}} = \frac{1}{{2\sqrt 2 {{\left( {1 - {\kappa _g} + {\kappa _g}^2} \right)}^{\frac{3}{2}}}}}\left( {{\lambda _1}\gamma  + {\lambda _2}t + {\lambda _3}d} \right)
\end{equation}
where

\[\begin{array}{l}
{\lambda _1} =  - {\kappa _g}^\prime  + 2{\kappa _g}{\kappa _g}^\prime  - 2 + 4{\kappa _g} - 4{\kappa _g}^2 + 2{\kappa _g}^3\\
{\lambda _2} =  - {\kappa _g}^\prime  - {\kappa _g}{\kappa _g}^\prime  - 2 - 4{\kappa _g}^2 + 2{\kappa _g} + 2{\kappa _g}^3 - 2{\kappa _g}^4\\
{\lambda _3} =  - {\kappa _g}{\kappa _g}^\prime  + 2{\kappa _g} - 4{\kappa _g}^2 + 2{\kappa _g}^\prime  + 4{\kappa _g}^3 - 2{\kappa _g}^4.
\end{array}\]
Substituting the equation (3.14) into equation (3.16), we reach

\begin{equation}\label{2.1}
{t_\beta }^\prime  = \frac{{\sqrt 3 }}{{4{{\left( {1 - {\kappa _g} + {\kappa _g}^2} \right)}^2}}}\left( {{\lambda _1}\gamma  + {\lambda _2}t + {\lambda _3}d} \right).
\end{equation}
Using the equations (3.13) and (3.15), we have

\begin{equation}\label{2.1}
{d_\beta } = \beta  \wedge {t_\beta } = \frac{1}{{\sqrt 6 \sqrt {1 - {\kappa _g} + {\kappa _g}^2} }}\left( {\left( {2{\kappa _g} - 1} \right)\gamma  + \left( { - 1 - {\kappa _g}} \right)t + \left( {2 - {\kappa _g}} \right)d} \right).\\
\end{equation}
From the equation (3.17) and (3.18), the geodesic curvature of  $\beta(s^{*})$ is

\[\begin{array}{l}
{\kappa _g}^\beta  = \left\langle {{t_\beta }^\prime ,{d_\beta }} \right\rangle \\\\
\,\,\,\,\,\,\,\,\,\,\,=\frac{1}{{4\sqrt 2 {{\left( {1 - {\kappa _g} + {\kappa _g}^2} \right)}^{\frac{3}{2}}}}}\left( {{\lambda _1}\left( {2{\kappa _g} - 1} \right) + {\lambda _2}\left( { - 1 - {\kappa _g}} \right) + {\lambda _3}\left( {2 - {\kappa _g}} \right)} \right).
\end{array}\]

\subsection{Example}
Let us consider the unit speed spherical curve:
\[\gamma \left( s \right) = \left\{ {cos\left( s \right){\rm{ }}tanh\left( s \right),{\rm{ }}sin\left( s \right){\rm{ }}tanh\left( s \right),{\rm{ }}sech\left( s \right)} \right\}.\]\\
It is rendered in Figure 1.

\begin{center}
\begin{figure}[htp]
\hfil\scalebox{1}{\includegraphics[width=0.5\textwidth]{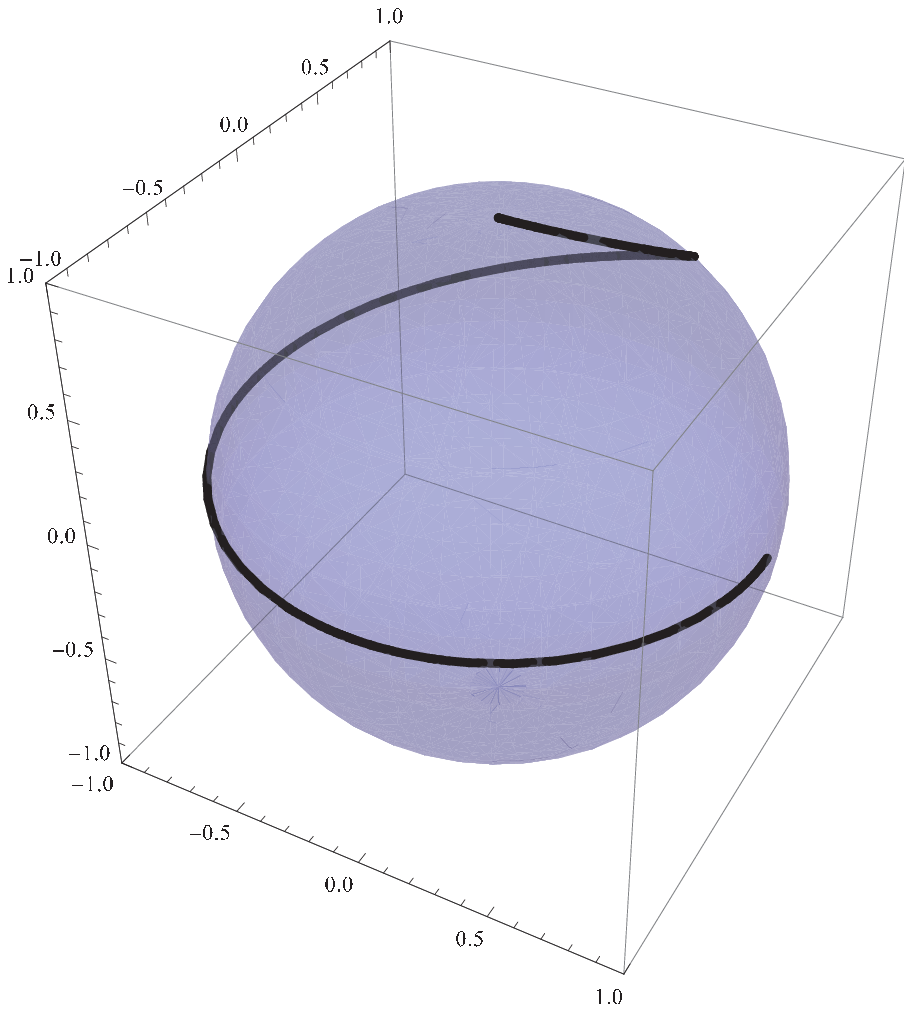}}\hfil
\caption{$\gamma=\gamma(s)$}
\end{figure}
\end{center}
In terms of definitions, we obtain Smarandache curves according to Sabban frame on $S^2$, see Figures 2 - 4.

\begin{center}
\begin{figure}[htp]
\hfil\scalebox{1}{\includegraphics[width=0.5\textwidth]{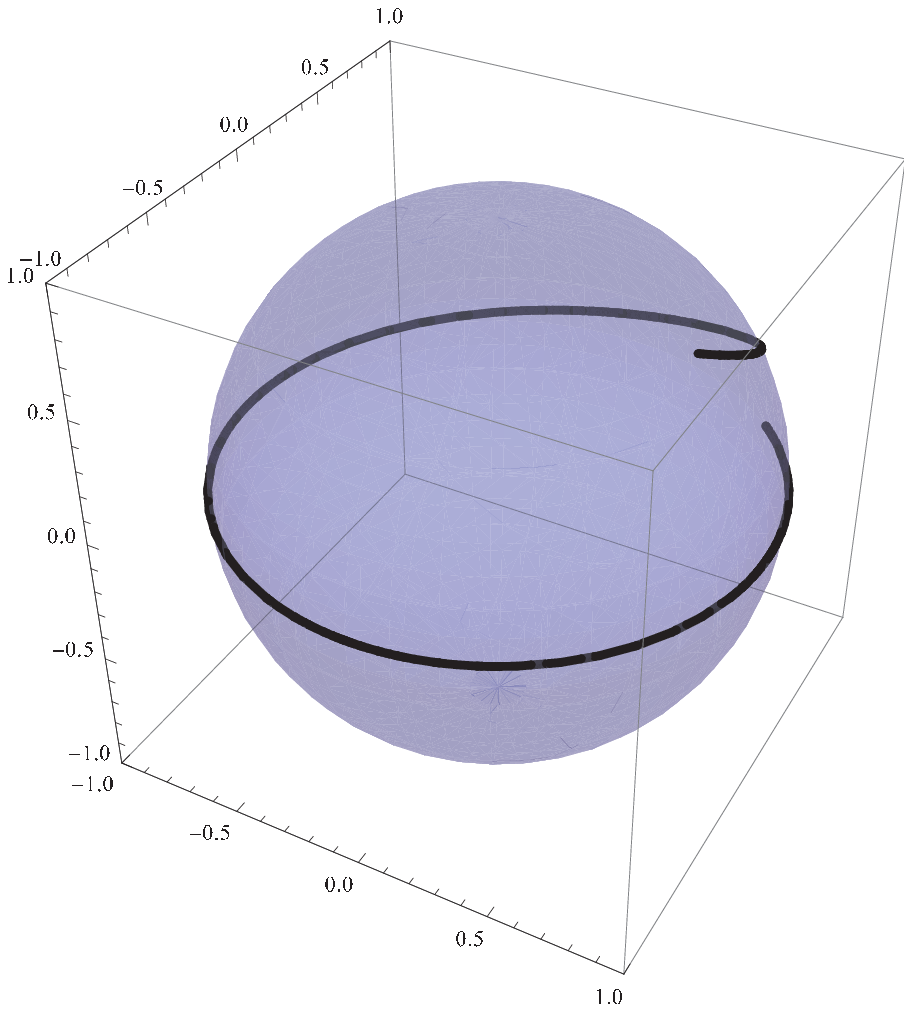}}\hfil
\caption{$\gamma t$ - Smarandache Curve}
\end{figure}
\end{center}

\begin{center}
\begin{figure}[htp]
\hfil\scalebox{1}{\includegraphics[width=0.5\textwidth]{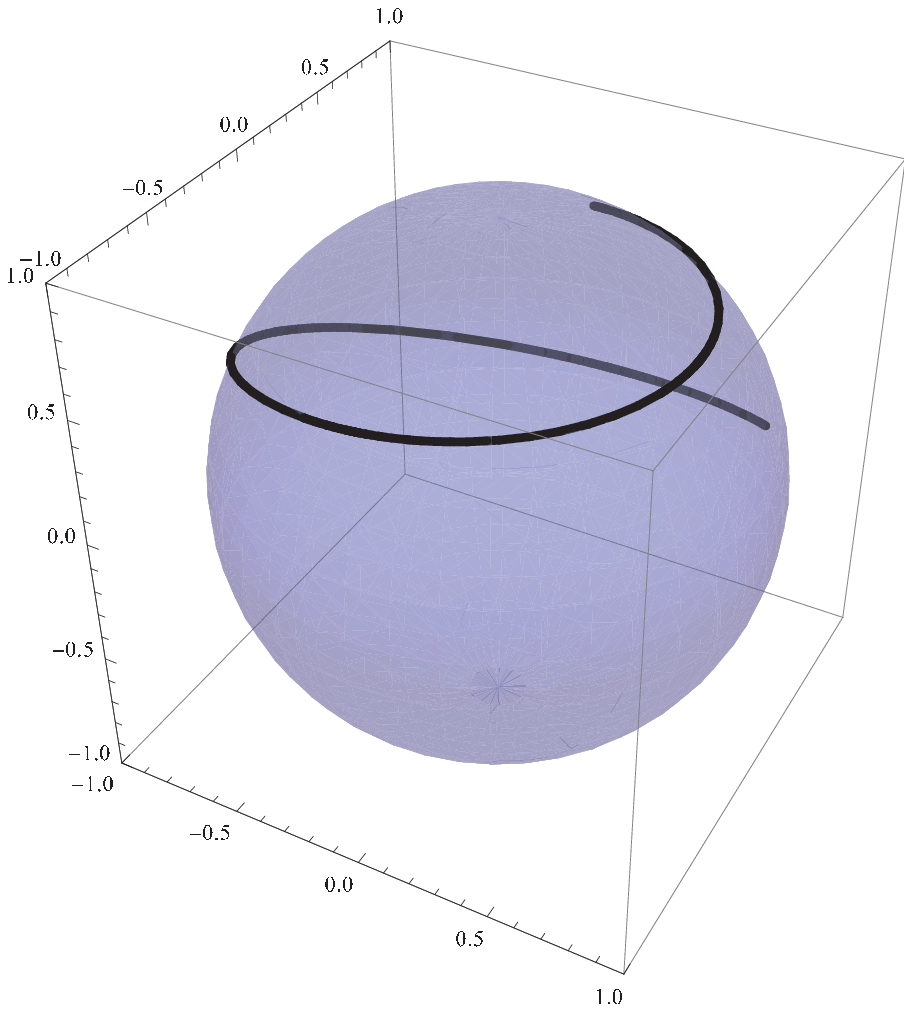}}\hfil
\caption{$td$ - Smarandache Curve}
\end{figure}
\end{center}

\begin{center}
\begin{figure}[htp]
\hfil\scalebox{1}{\includegraphics[width=0.5\textwidth]{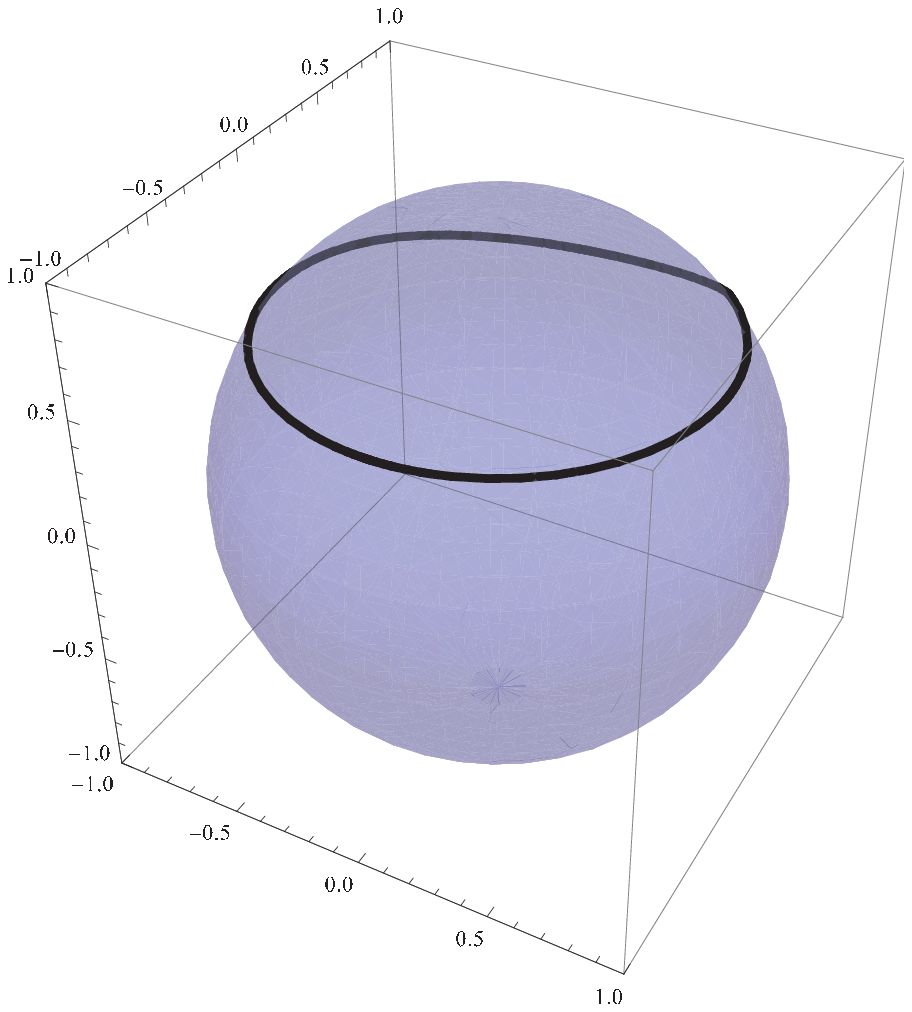}}\hfil
\caption{$\gamma td$ - Smarandache Curve}
\end{figure}
\end{center}

\end {document}